%&amstex          
\input amstex\documentstyle{amsppt}  
\pagewidth{12.5cm}\pageheight{19cm}\magnification\magstep1
\topmatter
\title{On bases of certain Grothendieck groups}\endtitle
\author G. Lusztig\endauthor
\address{Department of Mathematics, M.I.T., Cambridge, MA 02139}\endaddress
\thanks{Supported by NSF grant DMS-2153741}\endthanks
\endtopmatter   
\document
\define\lf{\lfloor}
\define\rf{\rfloor}

\define\mpb{\medpagebreak}

\define\sqc{\sqcup}

\define\bin{\binom}
\define\op{\oplus}
   
\define\part{\partial}
\define\emp{\emptyset}

\define\ra{\rangle}
\define\n{\notin}

\define\m{\mapsto}
\define\do{\dots}
\define\la{\langle}

\define\sub{\subset}    

\define\T{\times}
\define\ti{\tilde}
\define\nl{\newline}
\redefine\i{^{-1}}

\define\un{\underline}
\define\ov{\overline}

\define\supp{\text{\rm supp}}

\define\g{\gamma}

\define\e{\epsilon}
\define\et{\eta}
\define\io{\iota}

\redefine\t{\tau}

\redefine\aa{\bold a}

\define\NN{\bold N}

\define\SS{\bold S}

\define\ZZ{\bold Z}

\define\cc{\Cal C}

\define\ce{\Cal E}

\define\cj{\Cal J}

\define\cp{\Cal P}

\define\cs{\Cal S}

\define\cx{\Cal X}

\define\tE{\ti E}

\head Introduction\endhead
\subhead 0.1\endsubhead
Let $\cs$ be the set of (isomorphism classes of irreducible) unipotent representations of 

(a) a symplectic or odd special orthogonal group over a finite field, or

(b) an even split full orthogonal group over a finite field, or

(c) an even non-split full orthogonal group over a finite field.
\nl
(In cases (b),(c) we say that an irreducible
representation is unipotent if its restriction to the corresponding
special orthogonal group is unipotent; we further assume that this restriction is irreducible.)

Let $W$ be the corresponding Weyl group and let $Ce(W)$ be the set of two-sided cells of $W$,
which in cases (b),(c) are stable under the graph automorphism of $W$ induced by an element
in the full orthogonal group which is not in the special orthogonal group.

From \cite{L84} (or a slight extension) one deduces a natural partition

$\cs=\sqc_{c\in Ce(W)}\cs_c$.
\nl
Thus, the Grothendieck group of $\cs$ is a direct
sum $\op_{c\in Ce(W)}\ZZ[\cs_c]$ where for a finite set $Y$, $\ZZ[Y]$ denotes the free abelian
group with basis $Y$.

We now fix $c\in Ce(W)$. In
\cite{L20, L22} a new basis of $\ZZ[\cs_c]$ with strong positivity properties with respect to Fourier
transform was defined in case (a).
Here we shall call it the ``second basis''.
In case (b) there were two
versions of such a basis in \cite{L20} and the second one was
adopted in \cite{L22} under the name ``second basis''.

In this paper we give a somewhat different presentation
and refinements of the results of \cite{L20, L22} and extend them to include case (c).

From \cite{L84} (or a slight extension), to $c$ one can attach a pair of finite subsets $U'\sub U$ of
$\NN$ with $|U-U'|$ odd in case (a) and even in cases (b),(c), so that $\cs_c$ is identified with the
set $Sy(U',U)$ of ``symbols'', that is ordered pairs $(S,T)$ where $S,T$ are finite subsets of $U$
such that $S\cup T=U,S\cap T=U',0\n U'$ and

(d) $|S|-|T|$ is in $2\ZZ+1$ in case (a), in $4\ZZ$ in case (b), in $4\ZZ+2$ in case (c).
\nl
(The cardinal of a finite set $Y$ is denoted by $|Y|$.)
We can identify $Sy(U',U)$ with $Sy(\emp,[1,D+1])$ where
$D+1=|U-U'|$. (For $i,j$ in $\ZZ$ we set

$[i,j]=\{z\in\ZZ,i\le z\le j\}$.)
\nl
Namely, to $(S,T)\in Sy(U',U)$ we associate $(S',T')\in Sy(\emp,[1,D+1])$
where $S',T'$ are the images of $S-U',T-U'$ under the unique order preserving bijection
$U-U'@>>>[1,D+1]$. In this way $\ZZ[\cs_c]$
becomes $\ZZ[Sy(\emp,[1,D+1])]$ and our task becomes that of
defining a second basis of $\ZZ[Sy(\emp,[1,D+1])]$. 

We shall write $Sy_D$ (resp. $Sy^+_D,Sy^-_D$) instead of $Sy(\emp,[1,D+1])$ in case (a) (resp.
(b),(c)). Thus,

$Sy_D$ is the set of ordered pairs $(S,T)$ of disjoint subsets of $[1,D+1]$ such that
$S\cup T=[1,D+1]$, $|S|-|T|\in2\ZZ+1$, $D$ even $\ge0$;

$Sy_D^+$ is the set of ordered pairs $(S,T)$ of disjoint subsets of $[1,D+1]$ such that
$S\cup T=[1,D+1]$, $|S|-|T|\in4\ZZ$, $D$ odd $\ge1$;

$Sy_D^-$ is the set of ordered pairs $(S,T)$ of disjoint subsets of $[1,D+1]$ such that
$S\cup T=[1,D+1]$, $|S|-|T|\in4\ZZ+2$, $D$ odd $\ge1$.

Under our identification, the partition of $\cs_c$ according
to Harish-Chandra series corresponds to the partition

$Sy_D=\sqc_{s\in2\NN+1}Sy_D^s$, $Sy_D^s=\{(S,T)\in Sy_D;abs(|S|-|T|)=s\}$ in case (a);

$Sy^+_D=\sqc_{s\in4\ZZ}Sy_D^s$, $Sy_D^s=\{(S,T)\in Sy_D^+;|S|-|T|=s\}$ in case (b);

$Sy^-_D=\sqc_{s\in4\ZZ+2}Sy_D^s$, $Sy_D^s=\{(S,T)\in Sy_D^-;|S|-|T|=s\}$ in case (c).
\nl
Here, the absolute value of an integer $z$ is denoted by $abs(z)$.

\head 1. The second basis\endhead
\subhead 1.1\endsubhead
Let $F$ be the field consisting of two elements. 
Let $D\ge0$ be an integer. We set $N=D+1$ if $D$ is even, $N=D+2$ if $D$ is odd. Thus $N$ is odd.
When $N\ge3$, for any $k\in[1,D]$ we define an (injective) map $\io_k:[1,N-2]@>>>[1,N]$ by
$\io_k(i)=i$ for $i\in[1,k-1]$, $\io_k(i)=i+2$ for $i\in[k,N-2]$. Note that

$\io_k(1)<\io_k(2)<\do<\io_k(N-2)$ and $\io_k(i)=i\mod2$ for all $i\in[1,N-2]$.
\nl
Let $\tE_N$ be the set of subsets of $[1,N]$ viewed as an $F$-vector space in which the sum of $A,B$
is $(A\cup B)-(A\cap B)$. Let $E_N=\{X\in\tE_N;|X|=0\mod2\}$, a codimension $1$ subspace of $\tE_N$.
Let ${}^2E_N$ be the set of all $2$ element subsets of $E_N$.
When $N\ge3$, $k\in[1,D]$, we define an $F$-linear map $\tE_{N-2}@>>>\tE_N$ by
$\{i\}\m\{\io_k(i)\}$ for all $i\in[1,N-2]$; this map is denoted again by $\io_k$. It restricts
to an $F$-linear map $E_{N-2}@>>>E_N$.

A subset $\{i,j\}\in{}^2E_N$ will be often written as $ij$ if

$i<j$ and $i-j=1\mod2$ (we then say that $ij\in{}^2E'_N$) or if

$i>j$ and $i-j=0\mod2$ (we then say that $ij\in{}^2E''_N$).
\nl
In this way, ${}^2E_N={}^2E'_N\sqc{}^2E''_N$ is identified with a subset of $E_N\T E_N$.
For $ij\in{}^2E_N$ we set

$\lf i,j\rf=[i,j]$ if $i<j$, $i-j=1\mod2$ and

$\lf i,j\rf=[i,N]\sqc[1,j]$ if $i>j$, $i=j\mod2$.
\nl
We have

$|\lf i,j\rf|=j-i+1$ if $i<j$, $|\lf i,j\rf|=N-i+1+j$ if $i>j$.
\nl
Thus, $|\lf i,j\rf|=0\mod2$ that is, $\lf i,j\rf\in E_N$.

When $N\ge3$, $k\in[1,D]$ and $ij\in{}^2E_{N-2}$ we have, using the definitions, in $E_N$:

(a) $\lf\io_k(i),\io_k(j)\rf=\io_k(\lf i,j\rf)+c\{k,k+1\}$ where $c\in\{0,1\}$

\subhead 1.2\endsubhead
Let $\cp_N$ be the set of all unordered sequences $X_1,X_2,\do,X_t$ in ${}^2E_N$ such that
$X_a\cap X_b=\emp$ for any $a\ne b$. (We have necessarily $t\le(N-1)/2$.) 
For $B\in\cp_N$ let $\supp(B)=\cup_{ij\in B}\{i,j\}\sub[1,N]$ (this is a disjoint union) and
$B^1=B\cap{}^2E'_N$, $B^0=B\cap{}^2E''_N$; if $B^0\ne\emp$ we denote by $i_B$ the largest
number $i$ such that $ij\in{}^2E''_N$ for some $j$.

We define a subset $Pr_D$ of $\cp_N$ as follows: if $D=0$, $Pr_D$ consists of $Q_D^0=\emp$.
If $D$ is even $\ge2$, $Pr_D$ consists of $Q_D^0=\emp$ and of

$B=Q_D^t=\{\{D+1,1\},\{D,2\},\do,\{D+2-t,t\}\}$ with $t\in[2,D/2]$, $t$ even (we have $|B^0|=t$),

$B=Q_D^{-t}=\{\{D+1,1\},\{D,2\},\do,\{D+3-t,t-1\}\}$ with $t\in[2,(D+2)/2]$, $t$ even
(we have $|B^0|=t-1$).

If $D$ is odd, $Pr_D$ consists of $Q_D^{0,+}=\emp$, $Q_D^{0,-}=\{D+1,D+2\}$ and of  

$B=Q_D^{t,+}=\{\{D+1,2\},\{D,3\},\{D-1,4\},\do,\{D+3-t,t\}\}$ with $t$ even, $t\in[2,(D+1)/2]$
(we have $|B^0|=t-1$, $i_B\in2\NN$, $N\n\supp(B)$),

$B=Q_D^{-t,+}=\{\{D,1\},\{D-1,2\},\{D-2,3\},\do,\{D+2-t,t-1\}\}$ with $t$ even, $t\in[2,(D+1)/2]$
(we have $|B^0|=t-1$, $i_B\in2\NN+1$, $N\n\supp(B)$),

$B=Q_D^{-t,-}=\{\{D+2,1\},\{D+1,2\},\{D,3\},\do,\{D+4-t,t-1\}\}$ with $t$ even, $t\in[2,(D+3)/2]$,
(we have $|B^0|=t-1$, $i_B=N$),

$B=Q_D^{t,-}=\{\{D+1,D+2\},\{D,1\},\{D-1,2\},\do,\{D+1-t,t\}\}$ with $t$ even,
$t\in[2,(D-1)/2]$ (we have $|B^0|=t$, $i_B<N$, $N\in\supp(B)$).

For example,

$Pr_2$ consists of $\emp$ and $\{31\}$;

$Pr_4$ consists of $\emp$ and $\{51\},\{51,42\}$;

$Pr_6$ consists of $\emp$ and $\{71\},\{71,62\},\{71,62,53\}$;

$Pr_8$ consists of $\emp$ and $\{91\},\{91,82\},\{91,82,73\},\{91,82,73,64\}$;
                     
$Pr_1$ consists of $\emp$ and $\{31\},\{23\}$;

$Pr_3$ consists of $\emp$ and $\{42\},\{31\},\{51\},\{45\}$;

$Pr_5$ consists of $\emp$ and $\{62\},\{51\},\{71\},\{71,62,53\},\{67\},\{67,51,42\}$;

$Pr_7$ consists of $\emp$ and

$\{82\},\{82,73,64\},\{71\},\{71,62,53\},\{91\},\{91,82,73\},\{89\},\{89,71,62\}$.

When $N\ge3$ and $k\in[1,D]$ we define a map $I_k:\cp_{N-2}@>>>\cp_N$ by
$$(i_1j_1,i_2j_2,\do,i_tj_t)\m(\io_k(i_1)\io_k(j_1),\io_k(i_2)\io_k(j_2),\do,\io_k(i_t)\io_k(j_t),
\{k,k+1\}).$$

(This is well defined since $\io_k:[1,N-2]@>>>[1,N]$ is injective with image not containing $k,k+1$.) 

We define a subset $\cx_D$ of $\cp_N$ by induction on $D$ as follows.
We set $\cx_0=\cp_1$; it consists of $\emp$. We set $\cx_1=Pr_1\sqc\{12\}$. Assume now that $D\ge2$
so that $N\ge3$. Let $B\in\cp_N$. We say that $B\in\cx_D$ if either $B\in Pr_D$ or if there exists
$B'\in\cx_{D-2}$ and $k\in[1,D]$ such that $B=I_k(B')$.
We see that $Pr_D\sub\cx_D$ and that $I_k(\cx_{D-2})\sub\cx_D$ for $D\ge2,k\in[1,D]$.

In the remainder of this subsection we fix $B\in\cp_N$ and $\cj\sub[1,N]$ such that
either $\cj=[i,j]$ for some $i\le j$ or that $\cj=\emp$; define $e\in\{0,1\}$ by $|\cj|=e\mod2$.
We say that $\cj$ is {\it $e$-covered by $B$} if

(a) {\it there exists a sequence $a_1<b_1<a_2<b_2<\do<a_m<b_m$ in $[1,N]$ such that
$a_1b_1,a_2b_2,\do,a_mb_m$ are in $B^1$ and we have 

$\cj=[a_1,b_1]\sqc[a_2,b_2]\sqc\do\sqc[a_m,b_m]\sqc\cj_0$
\nl
where $\cj_0\sub\cj$ satisfies $|\cj_0|=e$.}

\subhead 1.3\endsubhead
For $B\in\cp_N$ we consider the following property:

$(P1)$ {\it There exists a sequence $i_1<i_2<\do<i_{2s}$ in $[1,N]$ such that
$B^0=\{i_{2s}i_1,i_{2s-1}i_2,\do,i_{s+1}i_s\}$ (it is automatically unique).}

\subhead 1.4 \endsubhead
For $B\in\cp_N$ we define $\et_B\in\{0,1\}$ as follows. If $D$ is even or if $D$ is odd and $|B^0|=0$
we have $|B^0|=\et_B\mod2$. If $D$ is odd and $|B^0|\ne0$ we set $\et_B=0$ if $N\in\supp(B)$ and
$i_B\ne N$ and $\et_B=1$ if $N\n\supp(B)$ or $i_B=N$. We consider the following property:

$(P2)$ {\it We have $|B^0|=\et_B\mod2$.}

\subhead 1.5\endsubhead
In this subsection we fix $B\in\cp_N$ such that $(P1)$ holds for $B$.
Let $i_1<\do<i_{2s}$ be as in $(P1)$. We consider the following property:

$(P3)$ {\it The following subsets of $[1,N]$ are $0$-covered by $B$:

$[i+1,j-1]$ for any $ij\in B^1$,

$[i_1+1,i_2-1],[i_2+1,i_3-1],\do,[i_{s-1}+1,i_s-1],[i_{s+1}+1,i_{s+2}-1],$

$[i_{s+2}+1,i_{s+3}-1],\do,[i_{2s-1}+1,i_{2s}-1]$ (if $s\ge1$).
\nl
(In particular any two consecutive terms of $i_1,i_2,\do,i_s$ have different parities, hence
any two consecutive terms of $i_{s+1},i_{s+2},\do,i_{2s}$ have different parities.) 
In addition, $[1,i_1-1],[i_{2s}+1,N]$ are $0$-covered by $B$ if $D$ is odd and $N\in\supp(B)$ or if
$D$ is even; if $D$ is odd and $N\n\supp(B)$, $i_B\in2\NN+1$, then $[1,i_1-1]$ is $0$-covered by $B$
and $[i_{2s}+1,N-1]$ is $1$-covered by $B$; if $D$ is odd and $N\n\supp(B)$, $i_B\in2\NN$, then
$[1,i_1-1]$ is $1$-covered by $B$ and $[i_{2s}+1,N-1]$ is $0$-covered by $B$ (again, if $s\ge1$).}

\subhead 1.6\endsubhead
Let $\ti\cx_D$ be the set of all $B\in\cp_N$ that satisfy $(P1),(P2),(P3)$. We show:

(a) {\it If $B\in\cx_D$ then $B\in\ti\cx_D$.}
\nl
We argue by induction on $D$. If $D\in\{0,1\}$ or if $D\ge2$ and $B\in Pr_D$, the result is obvious.
Thus we can assume that $D\ge2$ and $B=I_k(B')$ for some $k\in[1,D]$ and some $B'\in\cx_{D-2}$.
Then $(P1)$ for $B$ follows immediately from the analogous statement for $B'$.

To prove $(P2)$ for $B$ we can assume that $D$ is odd and $|B^0|\ne0$.   
Assume first that $|B^0|\in\{2,4,6,\do\}$. We have $|B^0|=|B'{}^0|$ so that
$|B'{}^0|\in\{2,4,6,\do\}$. From the induction hypothesis we see that $N-2\in\supp(B')$ and
$i_{B'}\ne N-2$. It follows that $N\in\supp(B)$ and $i_B\ne N$ as desired.
Next we assume that $|B^0|\in\{1,3,5,\do\}$. We have $|B^0|=|B'{}^0|$ so that
$|B'{}^0|\in\{1,3,5,\do\}$. From the induction hypothesis we see that $N-2\n\supp(B')$ or
$i_{B'}=N-2$. It follows that $N\n\supp(B)$ or $i_B=N$ as desired. We see that $(P2)$ holds for $B$.
It is easy to verify that if $(P3)$ holds for $B'$ then it holds for $B$.
This completes the proof of (a).

\subhead 1.7 \endsubhead
We show:

(a) {\it If $B\in\ti\cx_D$ then $B\in\cx_D$.}
\nl
We argue by induction on $D$. If $D\le1$, (a) is easily verified.
Now assume that $D\ge2$. Let ${}^*B$ be $B^0$ if $D$ is even and $B^0\cup\{N-1,N\}$ if $D$ is odd.
In the first part of the proof we assume that $B={}^*B$.

If $B^0=\emp$ then $B$ is either $\emp$ or $D$ is odd and $B=\{N-1,N\}$; in both cases we have
$B\in Pr_D$ and we are done. Thus we can assume that $B^0\ne\emp$. Let $i_1<\do<i_{2s}$ be as in
$(P1)$; note that $s\ge1$.
If $r\in\{1,2,\do,s-1,s+1,s+2,\do,2s\}$ and $i_{r+1}-i_r>1$, then by $(P3)$,
$[i_r+1,i_{r+1}-1]$ is $0$-covered by
$B$ and is nonempty, so that there exists $ab\in B^1$ such that $[a,b]\sub[i_r+1,i_{r+1}-1]$.
We have $ab\in B-{}^*B$, contradicting $B={}^*B$. We see that $i_{r+1}=i_r+1$ and
$$\align&(i_1,i_2,\do,i_s,i_{s+1},i_{s+2},\do,i_{2s})\\&
=(i_1,i_1+1,\do,i_1+s-1,i_{2s}-s+1,\do,i_{2s}-1,i_{2s}).\endalign$$
Assume now that $D$ is even. If $i_1\ge2$ then by $(P3)$, $[1,i_1-1]$ is $0$-covered by $B$
and is nonempty so that there exists $ab\in B^1$ with $[a,b]\sub[1,i_1-1]$.
This contradicts $B={}^*B$. We see that $i_1=1$.
If $N-i_{2s}\ge2$ then by $(P3)$, $[i_{2s}+1,N]$ is
$0$-covered by $B$ and is nonempty so that there exists $ab\in B^1$ with $[a,b]\sub[i_{2s}+1,N]$.
This contradicts $B={}^*B$. We see that $i_{2s}=N-1$ hence

$(i_1,i_2,\do,i_s,i_{s+1},i_{s+2},\do,i_{2s})=(1,2,\do,s,N-s,\do,N-2,N-1)$,
\nl
so that $B\in Pr_D$.

We now assume that $D$ is odd and $N\in\supp(B)$ or $N\n\supp(B),i_B\in2\NN+1$.
If $i_1\ge2$ then by $(P3)$, $[1,i_1-1]$ is $0$-covered by $B$ and is nonempty so that there exists
$ab\in B^1$ with $[a,b]\sub[1,i_1-1]$. We have $ab\in B-{}^*B$, contradicting $B={}^*B$. We see that
$i_1=1$.

We now assume that $D$ is odd and $N\n\supp(B)$, $i_B\in2\NN$.
If $i_1\ge3$ then by $(P3)$, $[1,i_1-1]$ is $1$-covered by $B$
and has at least $2$ elements, so that there exists $ab\in B^1$ with $[a,b]\sub[1,i_1-1]$.
We have $ab\in B-{}^*B$, contradicting $B={}^*B$. We see that $i_1\le2$.
Since $i_1=i_{2s}\mod2$ and $i_{2s}\in2\NN$ we see that $i_1=2$.

We now assume that $D$ is odd and $N\n\supp(B)$, $i_B\in2\NN$.
If $N-i_{2s}\ge3$, then by $(P3)$, $[i_{2s}+1,N-1]$ is $0$-covered by $B$ and is nonempty so that
there exists $ab\in B^1$ with $[a,b]\sub[i_{2s}+1,N-1]$. We have $ab\in B-{}^*B$, contradicting
$B={}^*B$. We see that $i_{2s}\ge N-2$ hence $i_{2s}=N-1$. Thus we have

$(i_1,i_2,\do,i_s,i_{s+1},i_{s+2},\do,i_{2s})=(2,\do,s+1,N-s,\do,N-2,N-1)$,
\nl
with $s$ odd (see $(P2)$), so that $B\in Pr_D$.

We now assume that $D$ is odd and $N\n\supp(B)$, $i_B\in2\NN+1$.
If $N-i_{2s}\ge3$, then by $(P3)$, $[i_{2s}+1,N-1]$ is $1$-covered by $B$
and has at least $2$ elements (since $i_{2s}+1=N-1\mod2$) so that there exists $ab\in B^1$
with $[a,b]\sub[i_{2s}+1,N-1]$.
We have $ab\in B-{}^*B$, contradicting $B={}^*B$. We see that $i_{2s}\ge N-2$;
since $i_{2s}=1\mod2, i_{2s}\ne N$, we must have $i_{2s}=N-2$ in this case. Thus we have

$(i_1,i_2,\do,i_s,i_{s+1},i_{s+2},\do,i_{2s})=(1,2,\do,s,N-s-1,\do,N-3,N-2)$,
\nl
with $s$ odd (see $(P2)$), so that $B\in Pr_D$.

We now assume that $D$ is odd and $N\in\supp(B)$, $i_B\ne N$.
If $N-i_{2s}\ge3$, then by $(P3)$, $[i_{2s}+1,N]$ is $0$-covered by $B$ and contains at least three
elements so that there exists $ab\in B^1$ other than $\{N-1,N\}$ with $[a,b]\sub[i_{2s}+1,N]$. 
We have $ab\in B-{}^*B$, contradicting $B={}^*B$. We see that $i_{2s}\ge N-2$.
As we have seen earlier, in this case we have $i_1=1$. Since $i_{2s}=i_1\mod2$ we see that
$i_{2s}$ is odd. Since $i_{2s}\ge N-2$ and $i_{2s}<N$ we see that $i_{2s}=N-2$.
Thus we have

$(i_1,i_2,\do,i_s,i_{s+1},i_{s+2},\do,i_{2s})=(1,2,\do,s,N-s-1,\do,N-3,N-2)$,
\nl
with $s$ even (see $(P2)$).
By $(P3)$, $[i_{2s}+1,N]=[N-1,N]$ is $0$-covered by $B$. It follows that $\{N-1,N\}\in B$.
We see that $B\in Pr_D$.

We now assume that $D$ is odd and $i_B=N$ hence $N\in\supp(B)$. In this case we have

$(i_1,i_2,\do,i_s,i_{s+1},i_{s+2},\do,i_{2s})=(1,2,\do,s,N-s+1,\do,N-1,N)$, 
\nl
with $s$ odd (see $(P2)$) so that $B\in Pr_D$.

We have thus proved that if $B={}^*B$ then $B\in Pr_D$; in particular we have $B\in\cx_D$.
We see that it is enough to prove (a) assuming that $B\ne{}^*B$.
Then we can find $ab\in B^1$ such that when $D$ is odd we have $\{a,b\}\ne\{N-1,N\}$.
We can assume in addition that $b-a$ is minimum possible. If $b-a>1$ then $b-a\ge3$ and
by $(P3)$, $[a+1,b-1]$ is $0$-covered by $B$ and is nonempty so that we can find $a'b'\in B^1$
with $[a',b']\sub[a+1,b-1]$; if $D$ is odd we have automatically $\{a',b'\}\ne\{N-1,N\}$.
This contradicts the minimality of $b-a$. We see that $b-a=1$, Thus there exists $k\in[1,D]$
such $\{k,k+1\}\in B$.
If $cd\in B-\{k,k+1\}$ then $\{c,d\}\cap\{k,k+1\}=\emp$ (since $B\in\cp_N$).
Hence there are unique $c',d'$ in $[1,N-2]$ such that $c=\io_k(c'),d=\io_k(d')$.
We have $c'd'\in{}^2E_{N-2}$. Let $B'\in\cp_{N-2}$ be the set consisting of all $c'd'$
associated as above to the various $cd\in B-\{k,k+1\}$. Note that $B=I_k(B')$. One can verify that
$(P1),(P2),(P3)$ hold for $B'$ since they hold for $B$. By the induction hypothesis we
have $B'\in\cx_{D-2}$. It follows that $B\in\cx_D$. This proves (a).

\subhead 1.8 \endsubhead
For $B\in\cx_D$ we show:

(a) {\it If $ab\in B,cd\in B$ are distinct (so that $\{a,b\}\cap\{c,d\}=\emp$) then
$\lf a,b\rf\cap\lf c,d\rf=\emp$ or $\lf a,b\rf\sub\lf c,d\rf$ or$\lf c,d\rf\sub\lf a,b\rf$.}
\nl
We argue by induction on $D$. If $D\in\{0,1\}$ or if $D\ge2$ and $B\in Pr_D$, the result is obvious.
Thus we can assume that $D\ge2$ and $B=I_k(B')$ for some $k\in[1,D]$ and some $B'\in\cx_{D-2}$.
Then $(a)$ for $B$ follows immediately from the analogous statement for $B'$.

\mpb

Let $B\in\cp_N$ and let $\cj\sub[1,N]$ with $\cj=[i,j]$ for some $i\le j$ be such that
$\cj$ is $e$-covered by $B$ (where $e=|\cj|\mod2,e\in\{0,1\}$) and let
$a_1<b_1<a_2<b_2<\do<a_m<b_m$ be as in 1.2(a). We show that:

(b) {\it the sequence $a_1<b_1<a_2<b_2<\do<a_m<b_m$ is unique.}
\nl
We argue by induction on $|\cj|$. If $|\cj|\le1$, the result is obvious. Now assume that
$|\cj|\ge2$ so that $\cj=[i,j]$ for some $i<j$.
Let $a'_1<b'_1<a'_2<b'_2<\do<a'_{m'}<b'_{m'}$ be a sequence with the same properties
as $a_1<b_1<a_2<b_2<\do<a_m<b_m$. Now $i$ equals $a_1$ or $a_1+1$ so that $i\in[a_1,b_1]$.
Similarly, $i\in[a'_1,b'_1]$. Using $(a)$ we see that $[a_1,b_1]=[a'_1,b'_1]$. Let
$\cj'=\cj-[a_1,b_1]$. We have $|\cj'|<|\cj|$.
The induction hypothesis is applicable to $\cj'$ instead of $\cj$; we see that (b) holds for $\cj$.

In the case where $e=1$, the unique element in $\cj$ which is not in
$[a_1,b_1]\sqc[a_2,b_2]\sqc\do\sqc[a_m,b_m]$ is said to be the {\it distinguished element} of $\cj$.

\subhead 1.9 \endsubhead
Assume  first that $D$ is even. We have a partition $\cx_D=\sqc_{t\in2\ZZ}\cx_D^t$ where
$\cx_D^t$ consists of all $B\in\cx_D$ such that $|B^0|=t$ (if $t\ge0$), $|B^0|=-t-1$ (if $t<0$).
The subsets $\cx_D^t$ are said to be the {\it pieces} of $\cx_D$. Note that when $t\in2\ZZ$,
$Q_D^t$ is the unique element in $Pr_D\cap\cx_D^t$.

In the remainder of this subsection we assume that $D$ is odd.
We have a partition $\cx_D=\cx_D^+\sqc\cx_D^-$ where

$\cx_D^+=\{B\in\cx_D;N\n\supp(B)\}$, $\cx_D^-=\{B\in\cx_D;N\in\supp(B)\}$.
\nl
For $t\in2\ZZ$ we define a subset $\cx_D^{t,+}$ of $\cx_D^+$ to be

$\{B\in\cx_D^+;|B^0|=0\}$ if $t=0$,

$\{B\in\cx_D^+;|B^0|=t-1,i_B\in2\NN\}$ if $t\in\{2,4,6,\do\}$,

$\{B\in\cx_D^+;|B^0|=-t-1,i_B\in2\NN+1\}$ if $t\in\{-2,-4,-6,\do\}$.
\nl
For $t\in2\ZZ$ we define a subset $\cx_D^{t,-}$ of $\cx_D^-$ to be

$\{B\in\cx_D^-;|B^0|=0\}$ if $t=0$,

$\{B\in\cx_D^-;|B^0|=t,i_B\ne N\}$ if $t\in\{2,4,6,\do\}$,

$\{B\in\cx_D^-;|B^0|=-t-1,i_B=N\}$ if $t\in\{-2,-4,-6,\do\}$.
\nl
Note that when $t\in2\ZZ$, $Q_D^{t,+}$ (resp. $Q_D^{t,-}$) is the unique element of
$Pr_D\cap\cx_D^{t,+}$ (resp. $Pr_D\cap\cx_D^{t,-}$).
From $(P3)$ we see that the subsets $\cx^{t,+},\cx^{t,-}$ of $\cx_D$ (with $t\in2\ZZ$) form a
partition of $\cx_D$; these subsets are said to be the {\it pieces} of $\cx_D$. 

\subhead 1.10 \endsubhead
In this subsection we asssume that $D$ is odd and that $B\in\cx_D$ is such that $|B^0|>0$ and
$N\n\supp(B)$. Let $i_1<\do<i_{2s}$ be as in (P1); here $s\ge1$.
If $i_{2s}\in2\NN+1$ (so that $i_{2s}+1=N-1\mod2$) then by $(P3)$, $[i_{2s}+1,N-1]$ is $1$-covered by
$B$; we denote by $u_B$ the distinguished element of $[i_{2s}+1,N-1]$. Note that
$[i_{2s}+1,u_B-1]$ and $[u_B+1,N-1]$ are $0$-covered by $B$.
If $i_{2s}\in2\NN$ (so that $i_1-1=1$) then by $(P3)$, $[1,i_1-1]$ is $1$-covered by $B$; we denote
by $u_B$ the distinguished element of $[1,i_1-1]$. Note that
 $[1,u_B-1]$ and $[u_B+1,i_1-1]$ are $0$-covered by $B$.

\subhead 1.11 \endsubhead
For any $B\in\cx_D$ we define $\un B\in E_N$ as follows. If $D$ is even we set $\un B=\emp$. If
$D$ is odd and $B\n\cup_{t\in2\ZZ-\{0\}}\cx_D^{t,+}$, we set $\un B=\emp$.
If $D$ is odd and $B\in\cx_D^{t,+}$ for some $t\in\{-2,-4,-6,\do\}$, we set
$\un B=\lf u_B,N\rf=[u_B,N]$; if $D$ is odd and $B\in\cx_D^{t,+}$ for some $t\in\{2,4,6,\do\}$,
we set $\un B=\lf N,u_B\rf=\{N\}\cup[1,u_B]$.

For any $B\in\cx_D$ we define $\e(B)\in E_N$ by

(a) $\e(B)=\sum_{ij\in B}\lf i,j\rf+\un B$.
\nl
(Sum in $E_N$.) We show:

(b) {\it If $N\ge3,k\in[1,D]$, $B'\in\cx_{D-2}$ then $\e(I_k(B'))=\io_k(\e(B'))+c\{k,k+1\}$
for some $c\in F$.}
\nl
We have:

$\e(I_k(B'))=\sum_{ij\in B'}\lf\io_k(i),\io_k(j)\rf+\{k,k+1\}+\un{I_k(B')}$,

$\io_k(\e(B'))=\sum_{ij\in B'}\io_k(\lf i,j\rf)+\io_k(\un{B'})$.
\nl
Using 1.1(a), we see that it is enough to show that

$\un{I_k(B')}=\io_k(\un{B'})+c_1\{k,k+1\}$ 
\nl
for some $c_1\in F$. If $D$ is even, both sides are zero (and $c_1=0$). Thus we can assume that $D$
is odd. From the definitions we have $u_B=\io_k(u_{B'})$, $N=\io_k(N-2)$ so that  the desired equality
follows again from 1.1(a). This proves (b).

\subhead 1.12\endsubhead
Define $\g:E_N@>>>\ZZ$ by $\g(X)=|X\cap(2\ZZ)|-|X\cap(2\ZZ+1)|$. Note that the image of
$\g$ is contained in $2\ZZ$. We show:

(a) {\it If $N\ge3,k\in[1,D]$, $B'\in\cx_{D-2}$ then $\g(\e(I_k(B')))=\g(\e(B'))$.}
\nl
Using 1.11(b), we see that it is enough to prove that for any $c\in F$ we have

$\g(\io_k(\e(B'))+c\{k,k+1\})=\g(\e(B'))$.
\nl
It is also enough to show that for any $X\in E_{N-2}$ we have $\g(\io_k(X)+c\{k,k+1\})=\g(X)$.
From the definitions we have $\g(\io_k(X))=\g(X)$ for any $X\in E_{N-2}$. It is then enough to show
that $\g(\io_k(X)+\{k,k+1\})=\g(X)$. From the definition we have $\io_k(X)\cap\{k,k+1\}=\emp$
hence $\g(\io_k(X)+\{k,k+1\})=\g(\io_k(X))+\g(\{k,k+1\})=\g(X)+\g(\{k,k+1\})=\g(X)$, since
$\g(\{k,k+1\})=0$. This proves (a).

\subhead 1.13\endsubhead
We now describe $\e(B)$ and $\g(\e(B)$ when $B\in Pr_D$. If $B=\emp$, then $\e(B)=\emp\in E_N$,
$\g(\e(B))=0$.

If $D$ is even and $B=Q_D^t$ with $t\in[2,D/2]$, $t$ even, we have

$\e(B)=\{2,4,6,\do,t,D+2-t,D+4-t,\do,D\}$, $\g(\e(B))=t$.

If $D$ is even and $B=Q_D^{-t}$ with $t\in[2,(D+2)/2]$, $t$ even, we have

$\e(B)=\{1,3,\do,t-1,D+3-t,D+5-t,\do,D+1\}$, $\g(\e(B))=-t$.

We now assume that $D$ is odd. If $B=Q_D^{0,-}$ then $\e(B)=\{D+1,D+2\}$ and $\g(\e(B))=0$.

If $B=Q_D^{t,+}$ with $t$ even, $t\in[2,(D+1)/2]$, we have

$\e(B)=\{2,4,6,\do,t,D+3-t,D+5-t,\do,D+1\}$, $\g(\e(B))=t$.

If $B=Q_D^{-t,+}$ with $t$ even, $t\in[2,(D+1)/2]$, then

$\e(B)=\{1,3,\do,t-1, D+2-t,D+4-t,\do,D\}$, $\g(\e(B))=-t$.

If $B=Q_D^{-t,-}$  with $t$ even, $t\in[2,(D+3)/2]$, then

$\e(B)=\{1,3,\do,t-1,D+4-t,D+6-t,\do,D+2\}$, $\g(\e(B))=-t$.

If $B=Q_D^{t,-}$ with $t$ even, $t\in[2,(D-1)/2]$, then

$\e(B)=\{2,4,6,\do,t,D+1-t,D+3-t,\do,D-1,D+1,D+2\}$, $\g(\e(B))=t$.

\subhead 1.14\endsubhead
We show:

(a) {\it If $N\ge3,k\in[1,D]$, $D$ odd, $B'\in\cx_{D-2}$, $B=I_k(B')\in\cx_D$, then
$N\in\e(B)$ if and only if $N-2\in\e(B')$.}
\nl
Recall from 1.11 that $\e(B)=\io_k(\e(B'))+c\{k,k+1\}$ for some $c\in\{0,1\}$. Note that
$N\n\{k,k+1\}$ so that we have $N\in\e(B)$ if and only if $N\in\io_k(\e(B'))$
and this happens if and only if $N-2\in\e(B')$. (From the definitions we see that for
$X\sub[1,N-2]$ we have $N\in\i_k(X)$ if and only if $N-2\in X$.) This proves (a).

\subhead 1.15\endsubhead
When $D$ is even we have a partition $E_N=\sqc_{t\in2\ZZ}\ce_D^t$ where

$\ce_D^t=\{X\in E_N;\g(X)=t\}$.
\nl
Now assume that $D$ is odd. We define a partition $E_N=\ce_D^+\sqc\ce_D^-$ by

$\ce_D^+=\{X\in E_N;N\n X\}$, $\ce_D^-=\{X\in E_N;N\in X\}$.
\nl
We have
$\ce_D^+=\sqc_{t\in2\ZZ}\ce_D^{t,+}$, $\ce_D^-=\sqc_{t\in2\ZZ}\ce_D^{t,}$ 
where

$\ce_D^{t,+}=\{X\in\ce_D^+;\g(X)=t\}$, $\ce_D^{t,-}=\{X\in\ce_D^-;\g(X)=t\}$.
\nl
We show:

(a) If $D$ is even and $B\in\cx_D^t$, $t\in2\ZZ$, then $\g(\e(B))=t$.

(b) If $D$ is odd and $B\in\cx_D^{t,+}$, $t\in2\ZZ$, then $\g(\e(B))=t$ and $N\n\e(B)$.

(c) If $D$ is odd and $B\in\cx_D^{t,-}$, $t\in2\ZZ$, then $\g(\e(B))=t$ and $N\in\e(B)$.
\nl
We argue by induction on $D$. For $D\in\{0,1\}$, the result is obvious.
When $B\in Pr_D$ this follows from 1.13. Assume now that $D\ge2$ and $B\n Pr_D$. We can find
$B'\in\cx_{D-2},k\in[1,D]$ such that $B=I_k(B')$. From the definitions we have
$B'\in\cx_{D-2}^t$ (in (a)), $B'\in\cx_{D-2}^{t,+}$ (in (b)), $B'\in\cx_{D-2}^{t,-}$ (in (c)).
By the induction hypothesis we have $\g(\e(B'))=t$ (in case (a),(b),(c)), $N\n\e(B')$
(in case (b)), $N\in\e(B')$ (in case (c)). Using 1.14(a) we deduce that $\g(\e(B))=t$ (in
case (a),(b),(c)). Using 1.11(b) we see that $N\n\e(B)$ in case (b) and
$N\in\e(B)$ in case (c). This completes the proof of (a),(b),(c).

\subhead 1.16\endsubhead
Let $t\in2\ZZ$. From 1.15(a),(b),(c) we see that $\e:\cx_D@>>>E_N$ restricts to a map

(a) $\cx_D^t@>>>\ce_D^t$
\nl
if $D$ is even, and to maps

(b) $\cx_D^{t,+}@>>>\ce_D^{t,+}$, $\cx_D^{t,-}@>>>\ce_D^{t,-}$
\nl
if $D$ is odd. Hence it restricts to maps

(c) $\cx_D^+@>>>\ce_D^+$, $\cx_D^-@>>>\ce_D^-$.

\subhead 1.17\endsubhead
For $B\in\cx_D$ let

(a) $\la B\ra$ be the subspace of $E_N$ spanned by the vectors $ij\in B$ (viewed as elements of $E_N$).
\nl
We show:

(b) $\e(B)\in\la B\ra$.
\nl
We argue by induction on $D$. If $D\in\{0,1\}$ the result is obvious. We now assume that $D\ge2$.
If $B\in Pr_D$ the result follows from 1.13. Thus we can assume that $B\n Pr_D$ so that we can find
$B'\in\cx_{D-2},k\in[1,D]$ such that $B=I_k(B')$. By 1.11(b) we have
$\e(B)=\io_k(\e(B'))+c\{k,k+1\}$ for some $c\in F$. From the definition of $I_k$ we see that
$\{k,k+1\}\sub\la B\ra$ and by the induction hypothesis we have $\e(B')\in\la B'\ra$. Thus it is
enough to prove that $\io_k(\la B'\ra)\sub\la B\ra$ or that for any $ij\in B'$ we have
$\io_k(\{i,j\})\in\la B\ra$; this follows from $\{\io_k(i),\io_k(j)\}\in B$.

\mpb

Let $\le_D$ be the transitive relation on $E_N$ generated by the relation for which
$X,X'$ in $E_N$ are related if $X\in\la\e\i(X')\ra$.

\proclaim{Theorem 1.18} (a) There is a unique bijection $\e':\cx_D@>>>E_N$ such that for any
$B\in\cx_D$ we have $\e'(B)\in\la B\ra$.

(b) We have $\e'(B)=\e(B)$ for any $B\in\cx_D$.

(c) The relation $\le_D$ is a partial order on $E_N$.

(d) The maps 1.16(a), 1.16(b), 1.16(c) are bijections.
\endproclaim
When $D$ is even, (a),(b),(c) can be deduced from the results of \cite{L20, L22}, see \S3.
The formula for $\e'$ given by (a) is simpler than the one in \cite{L20, L22}; the equivalence
of the two formulas is proved in \S3. The proof of (a),(b),(c)
for odd $D$ can be given along similar lines. Now (d) follows from (a),(b).

\subhead 1.19\endsubhead
For $X\in E_N$ we set $\cc(X)=[1,N]-X$ and

$X^*=(X\cap(2\ZZ+1))\cup(\cc(X)\cap(2\ZZ))\sub[1,N]$,

$X^{**}=(X\cap(2\ZZ))\cup(\cc(X)\cap(2\ZZ+1))=\cc(X^*)\sub[1,N]$.
\nl
We have

$|X^*|=|X\cap(2\ZZ+1)|+|(2\ZZ)\cap[1,N]|-|X\cap(2\ZZ)|=(N-1)/2-\g(X)$,
\nl
$|X^{**}|=(N+1)/2+\g(X)$. Hence $|X^{**}|-|X^*|=2\g(X)+1$.

Assume now that $D$ is even.
The assignment $X\m(X^*,X^{**})$ defines a bijection $E_N@>>>Sy_D$ (notation of 0.1);
it restricts for any $t\in2\ZZ$ to a bijection $\ce_D^t@>>>Sy_D^{abs(2t+1)}$.

For $X,X'$ in $E_N$ we set $M_{X,X'}=1$ if $X\in\la\e\i(X')\ra$ and $M_{X,X'}=0$, otherwise.
From 1.18 we see that $(M_{X,X'})$ is an upper triangular matrix with entries in $\{0,1\}$ and with
$1$ on diagonal. It follows that the elements

(a) $\sum_{X\in E_N}M_{X,X'}X'\in\ZZ[E_N]$ (for various $X'\in E_N$)
\nl
form a $\ZZ$-basis of $\ZZ[E_N]$, said to be the second basis. (This basis appears in \cite{L20} where it is called the new basis.)

Using the bijection $E_N@>>>Sy_D$ we see that the second basis of $\ZZ[E_N]$ can be viewed as a
$\ZZ$-basis of $\ZZ[Sy_D]$ which is also called the second basis.

We now assume that $D$ is odd.

If $X\in\ce_D^+$ then $(X^{**}-\{N\},X^*)\in Sy_D^+$ and
$|X^{**}-\{N\}|-|X^*|=2\g(X)$; the assignment $X\m(X^{**}-\{N\},X^*)$ defines
a bijection $\ce_D^+@>>>Sy_D^+$; it restricts for any $t\in2\ZZ$ to
a bijection $\ce_D^{t,+}@>>>Sy_D^{2t}$.

For $X,X'$ in $\ce_D^+$ we set $M^+_{X,X'}=1$ if $X\in\la\e\i(X')\ra$ and $M^+_{X,X'}=0$, otherwise.
From 1.18 we see that $(M^+_{X,X'})$ is an upper triangular matrix with entries
in $\{0,1\}$ and with $1$ on diagonal. It follows that the elements

(b) $\sum_{X\in\ce_D^+}M_{X,X'}X'\in\ZZ[\ce_D^+]$ (for various $X'\in\ce_D^+$)
\nl
form a $\ZZ$-basis of $\ZZ[\ce_D^+]$, said to be the second basis. (This basis appears in \cite{L22}.)
Using the bijection
$\ce_D^+@>>>Sy_D^+$, we see that the second basis of $\ZZ[\ce_D^+]$ can be viewed as a
$\ZZ$-basis of $\ZZ[Sy_D^+]$ which is also called the second basis.

If $X\in\ce_D^-$ then $(X^{**},X^*-\{N\})\in Sy_D^-$ and
$|X^{**}|-|X^*-\{N\}|=2\g(X)+2$; the assignment $X\m(X^{**},X^*-\{N\})$ defines a bijection
$\ce_D^-@>>>Sy_D^-$; it restricts for any $t\in2\ZZ$ to a bijection
$\ce_D^{t,-}@>>>Sy_D^{2t+2}$.

For $X,X'$ in $\ce_D^-$ we set $M^-_{X,X'}=1$ if $X\in\la\e\i(X')\ra$ and $M^-_{X,X'}=0$, otherwise.
From 1.18 we see that $(M^-_{X,X'})$ is an upper triangular matrix with entries
in $\{0,1\}$ and with $1$ on diagonal. It follows that the elements

(c) $\sum_{X\in\ce_D^-}M_{X,X'}X'\in\ZZ[\ce_D^-]$ (for various $X'\in\ce_D^-$)
\nl
form a $\ZZ$-basis of $\ZZ[\ce_D^-]$, said to be the second basis. 

Using the bijection $\ce_D^-@>>>Sy_D^-$, we see that the second basis of $\ZZ[\ce_D^-]$ can
be viewed as a $\ZZ$-basis of $\ZZ[Sy_D^-]$ which is also called the second basis.

\subhead 1.20\endsubhead
We have

$|Sy_D^s|=\bin{N}{(s+N)/2}$ if $D$ is even, $s\in2\NN+1$,

$|Sy_D^s|=\bin{N-1}{(s+N-1)/2}$, if $D$ is odd, $s\in2\ZZ$.

Here $\bin{a}{b}$ is defined to be $0$ if $b<0$ or if $b>a$. It follows that

$|\ce_D^t|=\bin{N}{(abs(2t+1)+N)/2}$ for $D$ even, $t\in2\ZZ$,

$|\ce_D^{t,+}|=\bin{N-1}{(2t+N-1)/2}$ for $D$ odd, $t\in2\ZZ$,

$|\ce_D^{t,-}|=\bin{N-1}{(2t+N+1)/2}$ for $D$ odd, $t\in2\ZZ$.

\head 2. Tables for $\cx_D$\endhead
\subhead 2.1\endsubhead
In this section we give tables describing $\cx_D$ and the map $\e:\cx_D@>>>E_N$ for $D=1,2,3,4,5,6,7$.
The table for $\cx_D$ is given by a list of elements of the various pieces of $\cx_D$; each such
element $B$ is written in the form $(?,?,\do,?,[?,?,\do,?])$ where each $?$ stands for an element of
$B$ and the $?$ inside the bracket $[,]$ are such that their sum in $E_N$ is equal to $\e(B)$.
The elements of $\cx_D$ are written in an order in which $B\in\cx_D$ preceeds $B'\in\cx_D$ whenever
$\e(B)\le_D\e(B')$.

\subhead 2.2. Table for $\cx_2$\endsubhead

$\cx_2^0$:

$([\emp]),([12]),([23])$

$\cx_2^{-2}$:

$([13]$

\subhead 2.3. Table for $\cx_4$\endsubhead

$\cx_4^0$:

$([\emp]), ([12]), ([23]), ([34]), ([45]),([12,34]), ([12,45])$,

$([23,45]), (23,[14]), (34,[25])$

\mpb

$\cx_4^{-2}$:

$([51]),([34,51]),([23,51]),(45,[31]),(12,[53])$

\mpb

$\cx_4^2$:

$(51,[42])$        

\subhead 2.4. Table for $\cx_6$\endsubhead

$\cx_6^0$:

$([\emp]), ([12]), ([23]), ([34]),([45]), ([56]), ([67]), ([12,45])$,

$([12,67]), ([23,45]), ([23,56]),([23,67]), ([34,67]), ([45,67]), ([12,34]), $

$([12,56]), ([34,56]), (23,[14]), (34,[25]), (45,[36]),$

$(56,[47]), ([12,34,56]), ([12,34,67]), (45,[12,36]), ([12,45,67]), ([23,45,67]),$

$(56,[12,47]),(56,[23,47]),(23,[56,14]), (23,[67,14]), (23,45,[16]), (25,[34,16]),$

$(34,[67,25]), (34,56,[27]),(36,[45,27])$

\mpb

$\cx_6^{-2}$:

$([71]), ([56,71]), ([45,71]), ([34,71]), ([23,71]), (67,[51]),$

$(12,[73]), ([23,56,71]),([34,56,71]),([23,45,71]),(67,[34,51]), $

$(12,[45,73]),(67,[23,51]),(12,[56,73]), (34,[25,71]), (45,[36,71]),$

$(45,67,[31]),(12,34,[75]),(12,67,[53]), (47,[56,31]), (14,[23,75])$

\mpb

$\cx_6^2$:

$(71,[62]),(71,[34,62]), (71,[45,62]), (71,23,[64]), (71,56,[42]),$

$(73,[12,64]), (51,[67,42])$

\mpb

$\cx_6^{-4}$:

$(62,[71,53])$       

\subhead 2.5. Table for $\cx_1$\endsubhead

$\cx_1^{0,+}$:

$([\emp]),([12])$

\mpb

$\cx_1^{0,-}$:

$([23])$   

\mpb

$\cx_1^{-2,-}$:

$[31]$

\subhead 2.6. Table for $\cx_3$\endsubhead

$\cx_3^{0,+}$:

$([\emp]),([12]),([23])$,

$([34]),([12,34]),(23,[14])$

\mpb

 $\cx_3^{-2,+}$:

$([31])$

\mpb

$\cx_3^{2,+}$:

$([42])$

\mpb

$\cx_3^{0,-}$:

$([45]),([12,45]),([23,45]),(34,[25])$

\mpb

$\cx_3^{-2,-}$:

$([51]),([51,34]),([51,23]),(12,[53])$      

\subhead 2.7. Table for $\cx_5$\endsubhead

$\cx_5^{0,+}$:

$([\emp]),([12]),([23]),([34]),([45]),([12,34])$,

$([12,45]),([23,45]),(23,[14]),(34,[25]),([56]),([12,56]),([34,56])$,

$([12,34,56]),([23,56]),(23,[14,56]),(45,[36]),(45,[12,36]),$

$(23,45,[16]),(25,[34,16])$

\mpb

$\cx_5^{-2,+}$:

$([51]),([34,51]),([23,51]),(45,[31]),(12,[53]),([31,56])$          
                        
\mpb

$\cx_5^{2,+}$:

$([62]),([34,62]),([45,62]),(56,[42]),(23,[64]),([12,64])$ 
         
\mpb

$\cx_5^{0,-}$:

$([67]),([12,67]),([23,67]),([34,67]),([45,67]),(56,[47])$,

$([12,34,67]),([12,45,67]),([23,45,67]),(23,[14,67]),(34,[25,67]),$

$(56,[12,47]),(56,[23,47]),(34,56,[27]),(36,[45,27])$

\mpb

$\cx_5^{-2,-}$:

$([71]),([23,71]),([34,71]),([45,71]),([56,71]),$

$([23,45,71]),([23,56,71]),([34,56,71]),(34,[25,71]),(45,[36,71]),$

$(12,[73]),(12,[45,73]),(12,[56,73]),(12,34,[75]),(14,[23,75])$            

\mpb

$\cx_5^{2,-}$:

$(51,[42,67])$            
                        
\mpb

$\cx_5^{-4,-}$:

$(62,[53,71])$

\subhead 2.8. Table for $\cx_7$\endsubhead

$\cx_7^{0,+}$:

$([\emp]),([12]),([23]),([34]),([45]),([56]),([67]),$

$([12,45]),([12,67]),([23,45]),([23,56]),([23,67]),([34,67]),([45,67]),$

$([12,34]),([12,56]),([34,56]),(23,[14]),(34,[25]),(45,[36]),$      

$(56,[47]),([12,34,56]),([12,34,67]),(45,[12,36]),([12,45,67]),$

$([23,45,67]),(56,[12,47]),(56,[23,47]),(23,[14,56]),(23,[14,67]),$

$(23,45,[16]),(25,[16,34]),(34,[25,67]),(34,56,[27]), (36,[27,45]),$   

$([78]),([12,78]),([23,78]),([34,78]),([45,78]),([56,78]),([12,34,78]),$

$([12,45,78]),([12,56,78]),([23,45,78]),([23,56,78]),([34,56,78]),$

$(67,[58]),(23,[14,78]),(34,[25,78]),(45,[36,78]),(67,[12,58]),$

$(67,[23,58]),(67,[34,58]),([12,34,56,78]),(23,[14,56,78]),(45,[12,36,78]),$

$(67,[12,34,58]),(45,67,[38]),(47,[38,56]),(23,45,[16,78]),(25,[16,34,78]),$

$(45,67,[12,38]),(23,67,[14,58],(23,45,67,[18]),(47,[12,38,56]),$

$(25,67,[18,34]),(23,47,[18,56]),(27,[18,34,56]),(27,45,[18,36])$  
                
\mpb

$\cx_7^{-2,+}$:

$([71]),([56,71]),([45,71]),([34,71]),([23,71]),$

$(67,[51]),(12,[73]),([23,71,56]),([34,71,56]), ([23,71,45]),(67,[34,51]),$

$(12,[73,45]),(67,[23,51]),(12,[73,56]),(34,[71,25]),(45,[71,36]), $

$(45,67,[31]),(12,34,[75]),(12,67,[53]),(47,[31,56]),(14,[75,23]),$

$([51,78]),([34,51,78]),(45,[31,78]),([23,51,78]),([31,56,78]),$

$(12,[53,78]),(67,[31,58])$      

\mpb

$\cx_7^{2,+}$:

$([82]),([82,34]),([82,45]),([82,56]),([82,67]),$

$([82,34,56]),([82,34,67]),([82,45,67]),(23,[84]),([84,12]),(23,[84,56]),$

$([84,12,56]),(23,[84,67]),([84,12,67]),(45,[82,36]),(23,45,[86]),$

$(45,[86,12]),(25,[86,34]),([86,12,34]),(23,[86,14]),(56,[82,47]),$

$(78,[62]),(78,[62,34]),(78,[62,45]),(78,56,[42]),(78,23,[64]),$

$(78,[64,12]),(58,[42,67])$   
             
 \mpb

$\cx_7^{-4,+}$:

$(62,[71,53])$   
             
\mpb

$\cx_7^{4,+}$:

$(73,[82,64])$ 
            
\mpb

$\cx_7^{0,-}$:

$([89]),([23,89]),([34,89]),([45,89]),([56,89]),$

$([67,89]),([78,69]),([12,34,89]),([12,45,89]),([12,56,89]),$

$([12,67,89]),([23,45,89]),([23,56,89]),([23,67,89]),([34,56,89]),$

$([34,67,89]),([45,67,89]),(23,[14,89]),(78,[12,69]),(78,[23,69]),$

$(78,[34,69]),(78,[45,69]),(34,[25,89]),(45,[36,89]),(56,[47,89]),$

$(56,78,[49]),(58,[67,49]),([12,34,56,89]),([12,34,67,89]),$

$([12,45,67,89]),([23,45,67,89]),(45,[12,36,89]),(78,[12,34,69]),$

$(78,[12,45,69]),(78,[23,45,69]),(56,[12,47,89]),(56,[23,47,89]),$

$(23,[14,56,89]),(23,[14,67,89]),(23,45,[16,89]),(34,56,[27,89]),$

$(23,78,[14,69]),(34,78,[25,69]),(25,[34,16,89]),(56,78,[12,49]),$

$(56,78,[23,49]),(36,[45,27,89]),(58,[12,67,49]),(58,[23,67,49]),$

$(34,56,78,[29]),(34,58,[67,29]),(36,78,[45,29]),(38,[45,67,29]),$

$(56,38,[47,29])$

\mpb

$\cx_7^{-2,-}$:

$([91]),([91,23]),([91,34]),([91,45]),([91,56]),$

$([91,67]),([91,78]),(12,[93]),([91,23,45]),([91,23,56]),([91,23,67]),$

$([91,23,78]),([91,34,56]),([91,34,67]),([91,34,78]),([91,45,67]),$

$([91,45,78]),([91,56,78]),([91,23,45,67]),([91,23,45,78]),([91,23,56,78]),$

$([91,34,56,78]),(34,[91,25]),(45,[91,36]),(56,[91,47]),(67,[91,58]),$

$(56,[91,23,47]),(67,[91,23,58]),(34,[91,25,67]),(34,[91,25,78]),$

$(67,[91,34,58]),(45,[91,36,78]),(34,56,[91,27]),(45,67,[91,38]),$

$(36,[91,27,45]),(47,[91,38,56]),(12,[93,45]),(12,[93,56]),(12,[93,67]),$

$(12,[93,78]),(12,[93,45,67]),(12,[93,45,78]),(12,[93,56,78]),$

$(12,56,[93,47]),(12,67,[93,58]),(12,34,[95]),(14,[23,95]),(12,34,[95,67]),$ 
              
$(12,34,[95,78]),(14,[23,95,67]),(14,[23,95,78]),(12,34,56,[97]),$

$(14,56,[23,97]),(12,36,[45,97]),(16,[23,45,97]),(16,34,[25,97])$

\mpb
                  
$\cx_7^{2,-}$:

$(71,[62,89]),(71,[62,89,45]),(71,[62,89,34]),$

$(71,23,[64,89]),(73,[64,89,12]),(71,56,[42,89]),(51,[42,67,89]),$

$(51,78,[42,69])$

\mpb

$\cx_7^{-4,-}$:

$(82,[73,91]),(82,[73,91,45]),(34,82,[75,91]),$

$(84,[75,91,23]),(12,84,[75,93]),(82,[73,91,56]), (67,82,[53,91]),$

$(62,[53,91,78])$

\head 3. Comparison with \cite{L22}\endhead
\subhead 3.1\endsubhead
In this section we assume that $D$ is even so that $N=D+1$.
The vectors $e_1=\{1,2\},e_2=\{2,3\},\do,e_D=\{D,D+1\}$ form a basis of
$E_N$. When $N\ge3$ we denote by $e'_1,e'_2,\do,e'_{D-2}$ the analogous vectors in $E_{N-2}$; for
$k\in[1,D],j\in[1,D-2]$ we have $\io_k(e'_j)=e_j$ if $j+1\le k-1$,
$\io_k(e'_j)=e_j+e_{j+1}+e_{j+2}$ if $j=k-1$, $\io_k(e'_j)=e_{j+2}$ if $k\le j,j+1\le N-2$.
Let $\SS^*_D$ be as in \cite{21, 3.1} where we take $\aa=N$. We have a unique bijection
$\t:\cx_D@>>>\SS^*_D$ which maps $ij$ to $\lf i,j\rf\in\SS^*_D$ for any $ij\in B$.
The bijection $\SS^*_D@>>>E_N$ in \cite{L22, 3.3(a)} becomes via $\t$ a bijection
$\ti\e:\cx_D@>>>E_N$ which, by \cite{L20, L22}, satisfies the requirement of 1.18(a) hence is equal to
$\e'$ in 1.18(a). We show that for any $B\in\cx_D$ we have

(a) $\ti\e(B)=\e(B)$,
\nl
so that 1.18(b) holds (for our $D$).
From \cite{L22, 3.2} we have
$$\ti\e(B)=\sum_{k\in[1,D]}[n_k]\{k,k+1\}+[n_{D+1}]\{1,N\}$$
(sum in $E_N$) where $n_k=|\{ij\in B;\{k,k+1\}\sub\lf i,j\rf\}|$ for $k\in[1,D]$ and
$n_{D+1}=|B^0|$ (for $n\in\NN$ we set $[n]=n(n+1)/2$); we view $[n_k]$ as integers $\mod2$.

From 1.11(a) we have
$$\e(B)=\sum_{k\in[1,N]}m_k\{k\}$$
(sum in $E_N$) where $m_k=|\{ij\in B;k\in\lf i,j\rf\}|$; we view $m_k$ as an integer $\mod2$.

For any $k\in[1,D+1]$ we set $k'=D+1,k''=2$ (if $k=1$),
$k'=k-1,k''=k+1$ (if $k\in[2,D]$), $k'=D,k''=1$ (if $k=D+1$). To prove (a) we must show that
$[n_{k}]+[n_{k'}]=m_k$ for $k\in[1,D+1]$ (equalities in $F$).

Let $Z_k=\{ij\in B;\{k',k,k''\}\sub\lf i,j\rf\}$,
$Z'_k=\{ij\in B;\{k',k\}\sub\lf i,j\rf,k''\n\lf i,j\rf\}$,
$Z''_k=\{ij\in B;\{k,k''\}\sub\lf i,j\rf,k'\n\lf i,j\rf\}$.
We have $n_k=|Z_k|+|Z''_k|,n_{k'}=|Z_k|+|Z'_k|,m_k=|Z_k|+|Z'_k|+|Z''_k|$.
From 1.8 we see that $|Z'_k|\in\{0,1\},|Z''_k|\in\{0,1\}$ and that at least one of
$|Z'_k|,|Z''_k|$ must be zero.
From 1.8 we see also that if $|Z_k|>0$ then there is a unique
$i_0j_0\in Z_k$ such that $\lf i_0,j_0\rf\sub\lf i,j\rf$ for any $ij\in Z_k$. Using $(P3)$ we see that

(b) if $k\in[2,D]$ then there exists $\ti i\ti j\in B^1$ such that $k\in[\ti i,\ti j]$ and
$[\ti i,\ti j]\subsetneqq\lf i_0,j_0\rf$;

(c) if $k=D+1$ then there exists $\ti i\ti j\in B^1$ such that $\ti j=k$
and $[\ti i,\ti j]\subsetneqq\lf i_0,j_0\rf$ (we use that $i_{2s}<D+1$, with notation of $(P1)$);

(d) if $k=1$ then there exists $\ti i\ti j\in B^1$ such that $\ti i=1$
and $[\ti i,\ti j]\subsetneqq\lf i_0,j_0\rf$ (we use that $i_1>1$, with notation of $(P1)$).

In case (b), by the minimality of $\lf i_0,j_0\rf$ we have $[\ti i,\ti j]\n Z_k$; since
$k\in[\ti i,\ti j]$ we must have either $[\ti i,\ti j]\in Z'_k$ or $[\ti i,\ti j]\in Z''_k$.
In case (c) we have $[\ti i,\ti j]\in Z'_k$; in case (d) we have $[\ti i,\ti j]\in Z''_k$.
Setting $A=|Z_k|,A'=|Z'_k|,A''=|Z''_k|$, we see that we have either

(i) $A>0,A'=1,A''=0$, or

(ii) $A>0,A'=0,A''=1$, or

(iii) $A=0,A'\in\{0,1\},A''\in\{0,1\}$.

The equality to be proved is

$(A+A')(A+A'+1)/2+(A+A'')(A+A''+1)/2=A+A'+A''$ (in $F$).
\nl
In case (i) this is the same as

$(A+1)(A+2)/2+A(A+1)/2=A+1$,
\nl
that is $A^2+2A+1=A+1$ (in $F$). This is
obvious. Now case (ii) is completely similar. In case (iii) we must
show that

$A'(A'+1)/2+A''(A''+1)/2=A'+A''$ (in $F$);
\nl
this is obvious when $A'\in\{0,1\},A''\in\{0,1\}$.

This completes the proof of (a).

\head 4. Even special orthogonal groups\endhead
\subhead 4.1\endsubhead
In this section we assume that we are in case 0.1(b) or 0.1(c). In this case, $\cs_c$ (see 0.1) admits
a fixed point free involution whose orbits are the pairs of unipotent representations in $\cs_c$ which
have isomorphic restrictions to the corresponding even special orthogonal group.
As in 0.1 we identify $\cs_c$ with $Sy_D^+$ (in case 0.1(b)) or with $Sy_D^-$ (in case 0.1(c)); here
$D$ is an odd integer $\ge1$. The
involution of $\cs_c$ becomes the fixed point free involution $(S,T)\m(T,S)$ of $Sy_D^+$ or $Sy_D^-$.
The set of orbits of this involution can be identified with the set $\cs'_c$ of unipotent
representations of the special orthogonal group attached to 0.1(b) or 0.1(c).

\subhead 4.2\endsubhead
Via the bijection $\ce_D^+@>>>Sy_D^+$ (resp. $\ce_D^-@>>>Sy_D^-$) in 1.19, the involution of
$Sy_D^+$ (resp. $Sy_D^-$) in 4.1 becomes the fixed point free involution $X\m X^!=X+[1,D+1]$ of
$\ce_D^+$, interchanging $\ce_D^{t,+},\ce_D^{-t,+}$ for any $t\in2\ZZ$, (resp. of $\ce_D^-$,
interchanging $\cx_D^{t,-},\cx_D^{-t-2,-}$ for any $t\in2\ZZ$).
Via the bijections 1.16(c) this becomes a fixed point free involution $B\m B^!$ of $\cx_D^+$,
interchanging $\cx_D^{t,+},\cx_D^{-t,+}$ for any $t\in2\ZZ$,
(resp. of $\cx_D^-$, interchanging $\cx_D^{t,-},\cx_D^{-t-2,-}$ for any $t\in2\ZZ$).

\subhead 4.3\endsubhead
We define a partition $\ce_D^{0,+}={}'\ce_D^{0,+}\sqc{}''\ce_D^{0,+}$ by
$${}'\ce_D^{0,+}=\{X\in\ce_D^{0,+};D+1\n X\},$$
$${}''\ce_D^{0,+}=\{X\in\ce_D^{0,+};D+1\in X\}.$$
Note that the involution $X\m X^!$ interchanges ${}'\ce_D^{0,+},{}''\ce_D^{0,+}$.

We define a partition $\cx_D^{0,+}={}'\cx_D^{0,+}\sqc{}''\cx_D^{0,+}$ by
$${}'\cx_D^{0,+}=\{B\in\cx_D^{0,+};D+1\n\supp(B)\},$$
$${}''\cx_D^{0,+}=\{B\in\cx_D^{0,+};D+1\in\supp(B)\}.$$
We show:

(a) {\it If $B\in{}'\cx_D^{0,+}$ then $\e(B)\in{}'\ce_D^{0,+}$.}
\nl
Indeed, since $\e(B)\in\la B\ra$ (see 1.17(b)), we see that $\e(B)$ is the sum in $E_N$ of certain
elements
$ij\in B$; now each such $ij$ satisfies $i\n N-1,j\n N-1$ so that $N-1\n\e(B)$, proving (a). We show:

(b) {\it If $B\in{}''\cx_D^{0,+}$ then $\e(B)\in{}''\ce_D^{0,+}$.}
\nl
By assumption there exists $ij\in B$ such that $i=N-1$ or $j=N-1$; the first possibility does not
occur since $B^0=\emp$ and $N\n\supp(B)$. Thus we have $\{h,N-1\}\in B$ for some $h\in[1,N-2]$;
since $B\in\cp_N$, such $h$ is in fact unique. If $ab\in B-\{h,N-1\}$ then we have $ab\in B^!$
since $B^0=\emp$ and $a<b<N-1$. We have $\e(B)=[h,N-1]+\sum_{ab\in B-\{h,N-1\}}[a,b]$ where the
last sum does not involve $N-1$; thus $N-1$ appears with coefficient $1$ in $\e(B)$ so that
$N-1\in\e(B)$, proving (b).

From  (a),(b) we see that the bijection $\cx_D^{0,+}@>>>\ce_D^{0,+}$ in 1.16(b) restricts to bijections
${}'\cx_D^{0,+}@>>>{}'\ce_D^{0,+}$, ${}''\cx_D^{0,+}@>>>{}''\ce_D^{0,+}$.
It follows that the involution $B\m B^!$ interchanges ${}'\cx_D^{0,+},{}''\cx_D^{0,+}$.

\subhead 4.4\endsubhead
One can verify that the following properties of $\le_D$ hold.

If $X\in\ce_D^{t,+},X'\in\ce_D^{t',+}$ are such that $X\le_DX'$ (with $t\in2\ZZ,t'\in2\ZZ$), then
we have $t=t'$ or $\max(t,-t)<\max(t',-t')$; if in addition $X'\in{}'\ce_D^{0,+}$, then
$X\in{}'\ce_D^{0,+}$.

If $X\in\ce_D^{t,-},X'\in\ce_D^{t',-}$ are such that $X\le_DX'$ (with $t\in2\ZZ,t'\in2\ZZ$), then
we have $t=t'$ or $\max(t,-t-2)<\max(t',-t'-2)$.

\subhead 4.5\endsubhead
Let

$\ce_D^{++}=\sqc_{t\in2\ZZ,t>0}\ce_D^{t,+}\sqc{}'\ce_D^{0,+}$,

$\ce_D^{+-}=\sqc_{t\in2\ZZ,t<0}\ce_D^{t,+}\sqc{}'\ce_D^{0,+}$,

$\ce_D^{-+}=\sqc_{t\in2\ZZ,t\ge0}\ce_D^{t,-}$,

$\ce_D^{--}=\sqc_{t\in2\ZZ,t<0}\ce_D^{t,-}$.

Note that each of $\ce_D^{++},\ce_D^{+-}$ is a set of representatives for the orbits of the
involution $X\m X^!$
of $\ce_D^+$ and that each of $\ce_D^{-+},\ce_D^{--}$ is a set of representatives for
the orbits of the involution $X\m X^!$ of $\ce_D^-$.

For $X,X'$ in $\ce_D^{++}$ let $M^{++}_{X,X'}=|\{Z\in\{X,X^!\};M^+_{Z,X'}=1\}|$.

For $X,X'$ in $\ce_D^{+-}$ let $M^{+-}_{X,X'}=|\{Z\in\{X,X^!\};M^+_{Z,X'}=1\}|$.

For $X,X'$ in $\ce_D^{-+}$ let $M^{-+}_{X,X'}=|\{Z\in\{X,X^!\};M^-_{Z,X'}=1\}|$.

For $X,X'$ in $\ce_D^{--}$ let $M^{--}_{X,X'}=|\{Z\in\{X,X^!\};M^-_{Z,X'}=1\}|$.

From 4.4 we see that $(M^{++}_{X,X'})$, $(M^{+-}_{X,X'})$, $(M^{-+}_{X,X'})$, $(M^{--}_{X,X'})$
are upper triangular matrices with entries in $\{0,1,2\}$ and with $1$ on diagonal.
It follows that

(a) $\sum_{X\in\ce_D^{++}}M^{++}_{X,X'}X'\in\ZZ[\ce_D^{++}]$ (for various $X'\in\ce_D^{++}$)
form a $\ZZ$-basis of $\ZZ[\ce_D^{++}]$ and

(b) $\sum_{X\in\ce_D^{+-}}M^{+-}_{X,X'}X'\in\ZZ[\ce_D^{+-}]$ (for various $X'\in\ce_D^{+-}$)
form a $\ZZ$-basis of $\ZZ[\ce_D^{+-}]$. 

It also follows that

(c) $\sum_{X\in\ce_D^{-+}}M^{-+}_{X,X'}X'\in\ZZ[\ce_D^{-+}]$ (for various $X'\in\ce_D^{-+}$)
form a $\ZZ$-basis of $\ZZ[\ce_D^{-+}]$ and

(d) $\sum_{X\in\ce_D^{--}}M^{--}_{X,X'}X'\in\ZZ[\ce_D^{--}]$ (for various $X'\in\ce_D^{--}$)
form a $\ZZ$-basis of $\ZZ[\ce_D^{--}]$.

\subhead 4.6\endsubhead
Let $\ov{\ce}_D^+$ (resp. $\ov{\ce}_D^-$) be the set of orbits of the fixed point free
involution $X\m X+[1,D+1]$ of $\ce_D^+$ (resp. $\ce_D^-$). The orbit maps
$\ce_D^+@>>>\ov{\ce}_D^+$, $\ce_D^-@>>>\ov{\ce}_D^-$ define bijections
$\ce_D^{++}@>>>\ov{\ce}_D^+$, $\ce_D^{+-}@>>>\ov{\ce}_D^+$, $\ce_D^{-+}@>>>\ov{\ce}_D^-$,
$\ce_D^{--}@>>>\ov{\ce}_D^-$ from which we get bijections
$\ce_D^{++}@>>>\cs'_c$, $\ce_D^{+-}@>>>\cs'_c$ (in case 0.1(b)) and
$\ce_D^{-+}@>>>\cs'_c$, $\ce_D^{--}@>>>\cs'_c$ (in case 0.1(c)).

Hence 4.5(a),(b) can be viewed as bases of the Grothendieck group $\ZZ[\cs'_c]$ (in case 0.1(b)) and
4.5(c),(d) can be viewed as bases of the Grothendieck group $\ZZ[\cs'_c]$ (in case 0.1(c)).

\enddocument